\documentclass[10pt,a4paper]{amsart}

\usepackage{dsfont}
\usepackage{amsmath}
\usepackage{xspace}
\usepackage{graphicx}
\usepackage{subfig}

\newcommand{\Gd}{\ensuremath{\mathds{G}_d}\xspace}
\newcommand{\Z}{\ensuremath{{\mathds{Z}}}\xspace}
\newcommand{\Zd}{\ensuremath{{\Z_d}}\xspace}
\newcommand{\C}{\ensuremath{\mathds{C}}\xspace}

\theoremstyle{plain}
\newtheorem{theorem}{Theorem}[section]
\newtheorem{lemma}[theorem]{Lemma}
\newtheorem{proposition}[theorem]{Proposition}
\newtheorem{corollary}[theorem]{Corollary}

\newtheorem{definition}[theorem]{Definition}

\newtheorem{remark}[theorem]{Remark}

\theoremstyle{definition}
\newtheorem{problem}[theorem]{Problem}

\numberwithin{equation}{section}

\newcommand{\nOffVs}[1]{\ensuremath{\kappa(#1)}}

\thanks{Both authors were
partially supported by UBACYT X042, CONICET PIP 5617 and ANPCyT PICT
20569, Argentina.}

\begin{document}

\author[A.~Dickenstein]{Alicia Dickenstein}
\email{alidick@dm.uba.ar}

\author[E.Tobis]{Enrique A. Tobis}
\email{etobis@dc.uba.ar}

\address{Departamento de Matem\'atica\\
FCEN, Universidad de Buenos Aires \\
(C1428EGA) Bue\-nos Aires, Argentina.}

\title{Additive edge labelings}

\begin{abstract}
  Let $G=(V,E)$ be a graph and $d$ a positive integer.  We study the
  following problem: for which labelings $f_E: E \to \Zd$ is there a
  labeling $f_V:V \to \Zd$ such that $f_E(i,j)\equiv f_V(i) + f_V(j)
  \pmod{d}$, for every edge $(i,j) \in E$? We also explore the
  connections of the equivalent multiplicative version to toric
  ideals. We derive a polynomial algorithm to answer these questions
  and to obtain all possible solutions.
\end{abstract}


\maketitle

\section{Introduction}
\label{sec:intro}

Graph labeling is a broad subject encompassing a myriad variants. In
its most general form, it involves assigning a value to each vertex or
each edge of a graph, subject to some restrictions. For an extensive
list of references on the subject, see the dynamic
survey~\cite{surveylabeling}.

A classic example of graph labeling is graph coloring.  Other
examples are harmonious labelings~\cite{harmonious} and felicitous
labelings~\cite{felicitous}. In the present work, we study a problem
similar to these last two, but dropping the one-to-one conditions and
allowing modular arithmetic over an arbitrary positive integer $d$.
The particular case $d=2$ is applied in~\cite{sontag} to the study of
monotone dynamical systems.

Let $G=(V,E)$ be a graph and let \Zd denote as usual the set of
integers modulo $d$. A function $f_E:E \to \Zd$ is called an
e-labeling of $G$ and a function $f_V:V \to \Zd$ is called a
v-labeling.  $(G,f_E)$ denotes the graph $G$ with its edges labeled
with $f_E$, and we say it is an \emph{e-labeled graph}.

In this paper, we answer completely the question of when a given 
labeling of the edges of $G$ with integers modulo $d$, admits a
labeling of the nodes of $G$ such that the label of each edge is the
sum, modulo $d$, of the labels of its vertices. More formally, we
study the following problem.

\begin{problem}
  \label{prob:additive}
  Let $(G,f_E)$ be an e-labeled graph. When is there a v-labeling
  $f_V$ of $G$ such that
  \begin{equation}
    f_E((v,v')) \equiv f_V(v) + f_V(v') \pmod{d}\label{eq:valid}
  \end{equation}
  holds for every edge $(v,v') \in E$? 
\end{problem}

\begin{definition}
  We say that an $f_V$ satisfying~(\ref{eq:valid}) is a \emph{valid
    v-labeling} of $(G,f_E)$. If such an $f_V$ exists, we say that
  $f_E$ is an \emph{additive e-labeling} of $G$. We also say that
  $(G,f_E)$ is an \emph{additively e-labeled graph}.
\end{definition}

Note that we are not imposing the restriction that adjacent vertices
have different labels.

\begin{figure}[ht]
  \centering
  \subfloat[]{%
    \label{sf:ejemploa}
  \includegraphics{additive1.eps}
  }
  \subfloat[]{
    \label{sf:ejemplob}
  \includegraphics{additive2.eps}
  }
  \subfloat[]{
    \label{sf:ejemploc}
  \includegraphics{additive3.eps}
  }
  \caption{\protect\subref{sf:ejemploa} An e-labeling of a graph with
    $\Z_3$ and \protect\subref{sf:ejemplob} a valid v-labeling of
    it.\protect\subref{sf:ejemploc} A non-additive e-labeling of a
    graph, again with $\Z_3$.}
\label{fig:example}
\end{figure} 

Once we know that an e-labeling $f_E$ of a graph $G$ is additive, we
can investigate how many valid v-labelings it admits. We denote this
number by \nOffVs{G,f_E}. For example, the graph of
Figure~\ref{sf:ejemploa}, with the edge labels in $\Z _3$, is additive
and admits a unique valid v-labeling over $\Z_3$, shown in
Figure~\ref{sf:ejemplob}, whereas the e-labeling of the graph of
Figure~\ref{sf:ejemploc} is not additive.

We characterize completely the existence of valid v-labelings in
Theorem~\ref{theo:compeqadd} and we compute \nOffVs{G, f_E} in
Theorem~\ref{theo:numbers}. In fact, we go beyond a theoretical
characterization. We present a polynomial algorithm to decide the
existence of valid v-labelings of an e-labeled graph $(G,f_E)$ in
Theorem~\ref{th:complexity}. We moreover show that we can enumerate
all such valid v-labelings in polynomial time. We reach our complexity
results in Section~\ref{sec:additivity} through the computation of the
Smith Normal Form (SNF)~\cite{polySmith} of the incidence matrix of
the graph~\cite{minorsIncidence} and our Theorem~\ref{theo:kernels}.

In Section~\ref{sec:additivity}, we comment on the equivalent
multiplicative version Problem~\ref{prob:multiplicative} of
Problem~\ref{prob:additive}, linking graphs and toric ideals. In
particular, we obtain in Theorem~\ref{theo:modulartoric} a modular
version of classic results on the implicitization of toric
parametrizations.

\section*{Acknowledgements} We are grateful to Eduardo Cattani for
calling our attention to the paper (\cite{sontag}) and suggesting us
the problem we study.  We are also indebted to Mar\'{\i}a
Ang\'elica Cueto for her useful comments and to Min Chih Lin for his
bibliographical help.

\section{Characterization of additive e-labelings}
\label{sec:char}

In this section we show that if a given e-labeling is additive, this
imposes restrictions on the cycles in $G$. Throughout this work, cycle
will not necessarily mean \emph{simple}
cycle. Theorem~\ref{theo:compeqadd} shows that these restrictions are
in fact sufficient.
 
If $C=(V,E)$ is a cycle of length $k$ in $G$, we number its nodes
``consecutively'' $v_1,\ldots,v_k$ and its edges $e_1,\ldots,e_k$,
where $e_i=(v_i,v_{i+1})$ for all $i < k$, and $e_k=(v_k,v_1)$.

\begin{definition}
  We say that an e-labeled graph $(G,f_E)$ has the \emph{even cycle
    property} if every cycle of even length in $G$, with edges
  $e_1,\ldots,e_{2k}$, satisfies
  \begin{equation} \label{eq:even}
  \sum_{l \text{ odd}}f_E(e_l) \equiv \sum_{l \text{ even}}f_E(e_l)
  \pmod{d}.
  \end{equation}
\end{definition}

\begin{definition}
  Let $d$ be an even positive integer. We say that an e-labeled graph
  $(G,f_E)$ has the \emph{odd cycle property} if every cycle of odd
  length in $G$, with edges $e_1,\ldots,e_{2k+1}$, satisfies
  \begin{equation}\label{eq:odd}
    \frac{d}{2} \sum_{l=1}^{2k+1}f_E(e_l) \equiv 0\pmod{d}.
  \end{equation}
\end{definition}
Equivalently, the odd cycle property holds if the sum
\[
\sum_{l=1}^{2k+1}f_E(e_l)
 \]
is an even number for all odd cycles in $G$.

Note that if $d=2$, then both properties mean that the number of $1$'s
in the edges of a cycle of any length is even. This case was studied
in~\cite{sontag} in a multiplicative setting as in
Section~\ref{sec:additivity}.

\begin{definition} \label{def:conditions}
  Let $(G,f_E)$ be an e-labeled graph. We say that $(G,f_E)$ is
  \emph{compatible} if the following conditions hold
  \begin{itemize}
  \item $d$ is odd and $(G,f_E)$ has the even cycle property.
  \item $d$ is even and $(G,f_E)$ has both the even and the odd cycle
    properties.
  \end{itemize}
\end{definition}

\begin{remark}
  The preceding definitions take into account \emph{all} the cycles of
  a graph, not just its simple cycles. The example in
  figure~\ref{fig:bowtie} shows two simple cycles joined by a
  vertex. Both cycles are e-labeled, and each of them is additive in
  isolation. However, they assign different labels to the shared
  vertex. This incompatibility only appears if we check non-simple
  cycles too.

  We show in Theorem~\ref{theo:kernels} that the number of
  conditions to be checked to ensure that $(G,f_E)$ is compatible is
  in fact ``small''.
  \begin{figure}[ht]
    \centering
    \includegraphics[scale=.8]{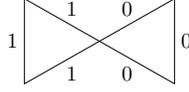}
    \caption{Two simple cycles joined by a vertex}
    \label{fig:bowtie}
  \end{figure}
\end{remark}

\begin{lemma}
  Let $(G,f_E)$, with $d$ even, be a connected e-labeled graph
  satisfying the even cycle property. Let $C$ be any odd cycle in $G$.
  Then $(G,f_E)$ satisfies the odd cycle property if and only if
  (\ref{eq:odd}) holds for $C$.
\end{lemma}
\begin{proof}
  We only need to prove one implication. Suppose that $(G,f_E)$
  satisfies the even cycle property, and that (\ref{eq:odd}) holds for
  $C$.  Let $C'$ be an odd cycle in $G$. Let $v \in C$ and $v' \in
  C'$. Since $G$ is connected, there is a path $P$ from $v$ to
  $v'$. Let $e_1,\ldots,e_{2k+1}$, $e'_1,\ldots,e'_{2s+1}$ and
  $e^P_1,\ldots,e^P_r$ be the edges of $C$, $C'$ and $P$, such that
  $v$ is a vertex of $e_1$ and of $e^P_1$, and such that $v'$ is a
  vertex of $e'_1$ and $e^P_r$. The even cycle property of $(G,f_E)$
  applied to the even cycle made up of $C$, $P$ from $v$ to $v'$, $C'$
  and then $P$ from $v'$ to $v$, implies that
  \begin{multline*}
    f_E(e_1) - f_E(e_2) + \cdots + f_E(e_{2k+1}) - f_E(e^P_1) +
    f_E(e^P_2) + \cdots + (-1)^rf_E(e^P_r) + \\ (-1)^{r+1}f_E(e'_1)+
    \cdots + (-1)^{r+2s+1}f_E(e'_{2s+1}) + (-1)^r f_E(e^P_r) + \cdots \ -
    f_E(e^P_1) \\ \equiv 0 \pmod{d}.
  \end{multline*}
  This is equivalent to
  \begin{multline*}
    f_E(e_1) - f_E(e_2) + \cdots + f_E(e_{2k+1}) - 2f_E(e^P_1) +
    2f_E(e^P_2) + \cdots + 
    \\ 2(-1)^rf_E(e^P_r) + (-1)^{r+1}f_E(e'_1) +
    \cdots + (-1)^{r+2s+1}f_E(e'_{2s+1}) \equiv 0 \pmod{d}.
  \end{multline*}
  If we multiply both sides by $d/2$, and since $d/2 \equiv -d/2
  \pmod{d}$, we get
  \begin{multline*}
    \frac{d}{2}\left(f_E(e_1) + f_E(e_2) + \cdots + f_E(e_{2k+1}) +
      f_E(e'_1) + \cdots + f_E(e'_{2s+1})\right) \equiv 0 \pmod{d}.
  \end{multline*}
  Using the odd cycle property of $(C,f_E)$, we get
  \begin{equation*}
    \frac{d}{2}\left(f_E(e'_1) +\cdots + f_E(e'_{2s+1})\right) \equiv 0
    \pmod{d},
  \end{equation*}
  which means that (\ref{eq:odd}) holds for $C'$ too. 
\end{proof}

We now show that compatibility is a necessary condition for
additivity.

\begin{lemma}
  \label{lem:evennecessary}
  If $(G,f_E)$ is an additive e-labeled graph, then $(G,f_E)$ has
  the even cycle property.
\end{lemma}
\begin{proof}
  Let $e_1,\ldots,e_{2k}$ be the edges of a cycle of even length in
  $G$. Recall that $e_i = (v_i,v_{i+1})$. Let $f_V$ be a v-labeling
  of $G$ satisfying~(\ref{eq:valid}). We have
  \begin{multline*}
    \sum_{l \text{ even}}f_E(e_l) \equiv \sum_{l \text{
        even}}(f_V(v_l) + f_V(v_{l+1})) \\
    \equiv \sum_{l \text{
        odd}}(f_V(v_l) + f_V(v_{l+1})) \equiv
    \sum_{l \text{ odd}}f_E(e_l)\pmod{d}.
  \end{multline*}
\end{proof}

\begin{lemma}
  \label{lem:oddnecessary}
  If $d$ is even, and $(G,f_E)$ is an additive e-labeled graph, then
  $G$ has the odd cycle property.
\end{lemma}
\begin{proof}
  Let $e_1,\ldots,e_{2k+1}$ be the edges of a cycle of odd length in
  $G$. Let $f_V$ be a v-labeling of $G$
  satisfying~(\ref{eq:valid}). We have
  \begin{equation*}
    \sum_{l=1}^{2k+1}\frac{d}{2}f_E(e_l) \equiv
    \sum_{l=1}^{2k+1}(\frac{d}{2}f_V(v_l) +
    \frac{d}{2}f_V(v_{l+1}))
    \equiv \sum_{l=1}^{2k+1}df_V(v_l) = 0\pmod{d}.
  \end{equation*}
\end{proof}

In fact, the compatibility conditions are sufficient for additivity.

\begin{theorem}
  \label{theo:compeqadd}
  An e-labeled graph $(G,f_E)$ is \emph{additive} if and only if it is
  \emph{compatible}.
\end{theorem}

Clearly, an e-labeled graph $(G,f_E)$ is additive if and only if every
connected component of G is additive with the labeling induced by
$f_E$. Also, to study the number of valid v-labelings of an e-labeled graph
$(G,f_E)$, we can assume that $G$ is a connected graph. Otherwise, if
$G_1, \ldots ,G_r$ are the connected components of $G$, we have
\begin{equation*}
  \nOffVs{G,f_E} = \prod_i\nOffVs{G_i,f_E}.
\end{equation*}

\begin{theorem}
  \label{theo:numbers}
  Let $(G,f_E)$ be a connected additive e-labeled graph.
  
  If $(G,f_E)$ has no odd simple cycles, $\nOffVs{G,f_E} = d$.
  
  If $(G,f_E)$ has at least one odd simple cycle, then
  \begin{itemize}
  \item if $d$ is odd then $\nOffVs{G,f_E} = 1$
  \item if $d$ is even then $\nOffVs{G,f_E} = 2$
  \end{itemize}
\end{theorem}

The proofs of these theorems occupy the next section.

\section{Proof of Theorems~\protect{\ref{theo:compeqadd}}
  and~\protect{\ref{theo:numbers}}}
\label{sec:proof}

Lemmas~\ref{lem:evennecessary} and~\ref{lem:oddnecessary} show that
compatibility is a necessary condition for additivity. We now turn our
attention to sufficient conditions and to the number of valid
v-labelings that an additive e-labeled graph admits, through a series
of preparatory lemmas.

\begin{lemma}
  \label{lem:determinedby1}
  Let $(G,f_E)$ be a connected additive e-labeled graph, and suppose
  that $f_V$, $f_{V'}$ are valid v-labelings of $(G,f_E)$. If there is
  $v \in V$ such that $ f_V(v) = f_{V'}(v)$, then $f_V=f_{V'}$.
\end{lemma}
\begin{proof}
  Let $v \in V$ be such that $f_V(v) = f_{V'}(v)$. We prove our lemma by
  induction on the distance from $v$. Let $v' \in V$ be at distance
  $0$ from $v$. Then, $v=v'$. Now, assume the lemma is true for all
  $v'$ at distance from $v$ smaller than $k$. Let $v'$ be at distance
  $k$. Let $\tilde{v} \in V$ be such that $d(v,\tilde{v}) = k-1$ and
  $d(\tilde{v},v') = 1$. Then, by our inductive hypothesis,
  $f_V(\tilde{v}) = f_{V'}(\tilde{v})$. Since $f_V$ and $f_{V'}$ are valid
  v-labelings, $f_V(v) = f_E(v,\tilde{v}) - f_V(\tilde{v}) =
  f_E(v,\tilde{v}) - f_{V'}(\tilde{v}) = f_{V'}(v)$.
\end{proof}

The previous lemma is important because it says that, given a
connected additive e-labeled graph, once we fix the label for one
vertex, the rest of the vertex labels are fixed. We use this result
later on. Furthermore, it shows that $\nOffVs{G,f_E} \leq d$.


\begin{definition}
  Given a simple cycle $C$ and three vertices $v$, $v'$ and $v''$ in
  $C$, we define $C[v,v',v'']$ as the simple path in $C$ from $v$ to
  $v'$ that contains $v''$. Conversely, $C[v,v',\overline{v''}]$ is
  the simple path from $v$ to $v'$ in $C$ that \emph{does not} contain
  $v''$(see Figure~\ref{fig:pathsincycle}.)
  \begin{figure}[ht]
    \centering
    \includegraphics{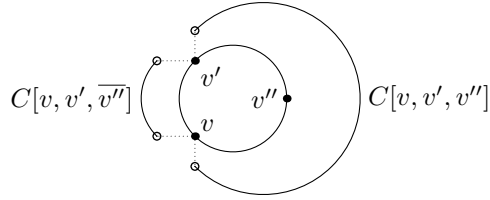}
    \caption{Two simple paths from $v$ to $v'$ in $C$.}
    \label{fig:pathsincycle}
  \end{figure}
\end{definition}

\begin{definition}
  Let $P$ be a 
  path with vertices $v_1, \ldots, v_k$ and edges
  $e_1=(v_1,v_2),\ldots,\\e_{k-1}=(v_{k-1},v_k)$. Let $f_E$ be an
  e-labeling of $P$. We define a function $\varphi_P: \Zd \to \Zd$.
  \begin{equation}
    \label{eq:varphi}
    \varphi_P(c) = (-1)^{k-1}c +
    \sum_{l=1}^{k-1}(-1)^{k-1-l}f_E(e_{l}) \pmod{d}.
  \end{equation}
  In other words, $\varphi_P(c)$ is the label that $v_k$ would have if we
  assigned label $c$ to $v_1$ and propagated it through $P$.
\end{definition}

\begin{remark}
  Let $(C,f_E)$ be an additive e-labeled simple cycle.  Let $v,v',v''$
  be in $C$ and set $C_1=C[v,v',v'']$ and
  $C_2=C[v,v',\overline{v''}]$. Let $f_V$ be a valid v-labeling of
  $(C,f_E)$. We have
  \begin{equation*}
    \varphi_{C_1}(f_V(v)) = \varphi_{C_2}(f_V(v)) = f_V(v').
  \end{equation*}
\end{remark}

We now prove Theorems~\ref{theo:compeqadd} and~\ref{theo:numbers} for
simple cycles of odd length.

\begin{lemma}
  \label{lem:oddcyclesadditive}
  If $(C,f_E)$ is a compatible e-labeled simple cycle of odd length
  then it is additive. If $d$ is odd, $\nOffVs{C,f_E} = 1$. If $d$ is
  even, $\nOffVs{C,f_E} = 2$.
\end{lemma}
\begin{proof}
  Let $v_1,\ldots,v_{2k+1}$ be the nodes of the cycle. Suppose that we
  have a valid v-labeling $f_V$. We want to see which are the possible
  values of $f_V(v_1)$. We need
  \begin{equation}
    \label{eq:conditionoddcycle}
    \varphi_{C[v_1,v_{2k+1},v_2]}(f_V(v_1)) \equiv
    \varphi_{C[v_1,v_{2k+1},\overline{v_2}]}(f_V(v_1)) \pmod{d}.
  \end{equation}
  We have
  \begin{equation*}
    \varphi_{C[v_1,v_{2k+1},v_2]}(f_V(v_1)) = (-1)^{2k}f_V(v_1) +
    \sum_{l=1}^{2k}(-1)^{2k-l}f_E(e_{l})\pmod{d},
  \end{equation*}
  and
  \begin{equation*}
    \varphi_{C[v_1,v_{2k+1},\overline{v_2}]}(f_V(v_1)) = f_E(e_{2k+1}) -
    f_V(v_1) \pmod{d}.
  \end{equation*}
  Merging these two expressions with~(\ref{eq:conditionoddcycle}) we get
  \begin{equation*}
    (-1)^{2k}f_V(v_1) + \sum_{l=1}^{2k}(-1)^{2k-l}f_E(e_{l}) \equiv
    f_E(e_{2k+1}) - f_V(v_1) \pmod{d}.
  \end{equation*}
  Since $2k$ is even, this expression is equivalent to
  \begin{equation}
    2f_V(v_1) \equiv \sum_{l=1}^{2k+1}(-1)^{l+1}f_E(e_l)
    \pmod{d}.\label{eq:condodd}
  \end{equation}

  \textbf{If $d$ is odd}, then $2$ is invertible modulo $d$ and
  equation~(\ref{eq:condodd}) has a unique solution. That implies that
  there is at most one possible value for $f_V(v_1)$. Since this value
  gives a valid v-labeling, there is a unique valid v-labeling of
  $(C,f_E)$.

  \textbf{If $d$ is even}, we use the odd cycle condition. Recall that
  this implies that the sum of the labels of the edges in the cycle is
  an even number.
   Since changing the sign of some summands does not alter the
  parity of a sum, the right side of~(\ref{eq:condodd}),
  \begin{equation*}
    \ell := \sum_{l=1}^{2k+1}(-1)^{l+1}f_E(e_l),
  \end{equation*}
  is also even.

  Equation~(\ref{eq:condodd}) is then of the form
  \begin{equation*}
    \label{eq:dioph2k}
    2X \equiv 2b \pmod{2c}.
  \end{equation*}
  This equation has exactly two solutions: $X = b$ and $X= b +
  c$. This means that $f_V(v_1)$ is either $\ell/2$ or $(\ell +
  d)/2$. Since these two values for $f_V(v_1)$ give valid v-labelings,
  our proof is complete.
\end{proof}

This proof allows us to deduce the following

\begin{corollary}
  \label{cor:difference}
  Let $(C,f_E)$ be an additive e-labeled simple cycle of odd length,
  with $d$ even. If $f_V$ and $f_{V'}$ are its two different valid
  v-labelings, then $f_V(v) \equiv f_{V'}(v) + d/2 \pmod{d}$ for all
  $v\in V$.
\end{corollary}

Let $(G,f_E)$ be an e-labeled graph. In the following proofs, we abuse
our notation. If $C$ is a subgraph of $G$, then $(C,f_E)$ stands for
the graph $C$ labeled with the restriction of $f_E$ to the edges of
$C$.

\begin{lemma}
  \label{lem:sumaround}
  Let $(G,f_E)$ be a compatible e-labeled connected graph. Let $C$ and
  $C'$ be two cycles of odd length in $G$. Let $e_1,\ldots,e_r$ be the
  edges of $C$ and $e'_1,\ldots,e'_s$ be the edges of $C'$. Assume
  that $C$ and $C'$ share at least one vertex $v_1$, such that both
  $e_1$ and $e'_1$ are incident to $v_1$. Then
  \begin{equation}
    \label{eq:sumOfTwoOddCycles}
    \sum_{l=1}^{r}(-1)^{r-l}f_E(e_l) \equiv
    \sum_{l=1}^{s}(-1)^{s-l}f_E(e'_l) \pmod{d}.
  \end{equation}
\end{lemma}
\begin{proof}
  Consider the cycle obtained by traversing
  $e_1,\ldots,e_r,e'_1,\ldots,e'_s$. Since $r$ and $s$ are odd, this
  cycle has even length. The compatibility hypothesis implies that
  \begin{equation}
    f_E(e_1) - f_E(e_2) + \cdots + f_E(e_r) - f_E(e'_1) + f_E(e'_2) -
    \cdots - f_E(e'_s) \equiv 0 \pmod{d}.
  \end{equation}
  But this means
  \begin{equation}
    \sum_{l = 1}^{r}(-1)^{r-l}f_E(e_l) - \sum_{l =
      1}^{s}(-1)^{s-l}f_E(e'_l) \equiv 0 \pmod{d},
  \end{equation}
  which is what we wanted to prove.
\end{proof}

\begin{proof}[Proof \textbf{of Theorems~\ref{theo:compeqadd}
    and~\ref{theo:numbers}}]
  Let $(G,f_E)$ be a compatible e-labeled graph. Without loss of
  generality, we can assume that it is connected. We prove the
  theorems by constructing a valid v-labeling of it.

  If $G$ has odd simple cycles, call one of them $C$. Choose a valid
  v-labeling $f$ of $(C,f_E)$. Pick a vertex $v$ in $C$ and set $\ell
  = f(v)$. If $G$ has no odd cycles, choose any vertex $v$ in $G$ and
  label it with any $\ell$ in $\Zd$.

  We build a valid v-labeling $f_V$ of $(G,f_E)$ by propagating the
  label of $v$ to the rest of the graph. For that, set $f_V(v) =
  \ell$. For any vertex $v' \in V$, choose a 
  path $P$ from $v$
  to $v'$ and set
  \[
  f_{V}(v') = \varphi_{P}(\ell),
  \]
  where $\varphi_P$ is as in~(\ref{eq:varphi}).
  We have to prove that $f_V$ is well defined and that it is a valid
  v-labeling of $(G,f_E)$.

  Given $v'$ and two simple paths $P_1$ and $P_2$ from $v$ to $v'$, we
  have to prove that
  \[
  \varphi_{P_1}(\ell) = \varphi_{P_2}(\ell).
  \]
  Let $e_1,\ldots,e_r$ and $e'_1,\ldots,e'_s$ be the edges of $P_1$
  and $P_2$, respectively, and assume that $v$ is an endpoint of $e_1$
  and $e'_1$. We call $C'$ the cycle formed by the union of $P_1$ and
  $P_2$.

  \textbf{If the sum of the lengths of $P_1$ and $P_2$ is even}, we
  can use the even cycle property of $(G,f_E)$ applied to $C'$. That
  is,
  \begin{equation*}
    f_E(e_1) - f_E(e_2) + \cdots + (-1)^{r+1}f_E(e_r) +
    (-1)^{r+2}f_E(e'_s) + \cdots - f_E(e'_1) \equiv 0 \pmod{d}.
  \end{equation*}
  This condition is equivalent to the identity
  \begin{equation}
    \label{eq:l}
    \sum_{l = 1}^r(-1)^{l}f_E(e_l) \equiv\sum_{l=1}^s(-1)^{l}f_E(e'_l) 
    \pmod{d}.
  \end{equation}
  We have that
  \begin{equation}
    \varphi_{P_1}(\ell) = (-1)^{r}\ell +
    \sum_{l=1}^{r}(-1)^{r-l}f_E(e_{l})\pmod{d},
  \end{equation}
  and
  \begin{equation}
    \varphi_{P_2}(\ell) = (-1)^{s}\ell +
    \sum_{l=1}^{s}(-1)^{s-l}f_E(e_{l})\pmod{d}.
  \end{equation}
  We must prove that $\varphi_{P_1}(\ell) =
  \varphi_{P_2}(\ell)$. Since the parity of $r$ and $s$ are the same,
  $(-1)^{s}\ell = (-1)^{r}\ell$, and we just need to prove that
  \begin{equation}
    \label{eq:pathsequal}
    \sum_{l=1}^{s}(-1)^{s-l}f_E(e_{l}) \equiv
    \sum_{l=1}^{r}(-1)^{r-l}f_E(e_{l}) \pmod{d}.
  \end{equation}
  If $r$ and $s$ are even, $(-1)^{r-l} = (-1)^{s-l} = (-1)^l$ for any
  integer $l$. Therefore,~(\ref{eq:l}) shows
  that~(\ref{eq:pathsequal}) holds. If $r$ and $s$ are odd,
  $(-1)^{r-l} = (-1)^{s-l} = (-1)^{l+1}$ for any integer $l$, and
  again~(\ref{eq:l}), this time multiplied by $-1$, shows
  that~(\ref{eq:pathsequal}) holds.

  \textbf{If $r$ is odd and $s$ is even}, the cycle $C'$ 
  has odd length. We need to prove that $\varphi_{P_1}(\ell)
  = \varphi_{P_2}(\ell)$, which is equivalent to
  \begin{equation}
    -\ell +
    \sum_{l=1}^{r}(-1)^{l+1}f_E(e_{l}) \equiv \ell +
    \sum_{l=1}^{s}(-1)^{l}f_E(e_{l}) \pmod{d}.
  \end{equation}
  This is the same as proving that
  \begin{equation}
    \label{eq:conditionodd}
    2\ell \equiv \sum_{l=1}^{r}(-1)^{l+1}f_E(e_{l}) +
    \sum_{l=1}^{s}(-1)^{l+1}f_E(e_{l}) \pmod{d}.
  \end{equation}
  The right side of~(\ref{eq:conditionodd}) is the alternating sum of
  the labels of the edges of the odd cycle $C'$, starting at $v$. By
  Lemma~\ref{lem:sumaround}, this sum is equal, modulo $d$, to the
  alternating sum of the labels of the edges of $C$, starting at
  $v$. By Lemma~\ref{lem:oddcyclesadditive}, this sum is equivalent to
  $2\ell$, which is what we needed to prove.
  
  We now know that $f_V$ is a well-defined labeling. We must show that
  it is also a valid v-labeling of $(G,f_E)$. That is, for each edge
  $(v',v'')$, 
  \begin{equation}
    \label{eq:additive}
    f_E((v',v'')) \equiv f_V(v') + f_V(v'') \pmod{d}.
  \end{equation}
  All the edges incident to $v$ satisfy~(\ref{eq:additive}) by the
  previous argument. Let $v'$ and $v''$ be two adjacent vertices in
  $G$, both different from $v$. Let $e$ be the edge between $v'$ and
  $v''$. Let $P_1$ and $P_2$ be paths from $v$ to $v'$ and $v''$,
  respectively. Let $e_1,\ldots,e_r$ and $e'_1,\ldots,e'_s$ be the
  edges of $P_1$ and $P_2$, respectively. We must prove that
  \begin{equation}
    \label{eq:additiveedge}
    \varphi_{P_1}(\ell) + \varphi_{P_2}(\ell) \equiv f_E(e) \pmod{d}.
  \end{equation}
  Consider the path $P'_2 = P_2 \cup \{e\}$. $P'_2$ and $P_1$ are two
  paths from $v$ to $v'$. If we write $e'_{s+1} = e$, we have just
  proved that
  \begin{equation}
    \label{eq:v'bothways}
    (-1)^r\ell + \sum_{l=1}^r(-1)^{r-l}f_E(e_l) \equiv 
    (-1)^{s+1}\ell + \sum_{l=1}^{s+1}(-1)^{s+1-l}f_E(e'_l) \pmod{d}.
  \end{equation}
  But the right side of~(\ref{eq:v'bothways}) can be split
  \begin{equation}
    \label{eq:lastedge}
    \sum_{l=1}^{s+1}(-1)^{s+1-l}f_E(e'_l) = f_E(e) +
    \sum_{l=1}^{s}(-1)^{s+1-l}f_E(e'_l).
  \end{equation}
  So joining~(\ref{eq:v'bothways}) and~(\ref{eq:lastedge}), we get
  \begin{equation}
    (-1)^r\ell + \sum_{l=1}^r(-1)^{r-l}f_E(e_l) + 
    (-1)^{s}\ell + \sum_{l=1}^{s}(-1)^{s-l}f_E(e'_l) \equiv  f_E(e)\pmod{d},
  \end{equation}
  which proves~(\ref{eq:additiveedge}).
\end{proof}

\section{Other results on Compatibility}

\begin{lemma}
  \label{lem:extension}
  Given a compatible e-labeled graph $(G,f_E)$, if we add any edge $e$
  to $G$, there is an extension of $f_E$ that assigns a label to $e$
  such that the resulting e-labeled graph is compatible.
\end{lemma}
\begin{proof}
  If we add an edge to a graph $G$, we can have three mutually
  exclusive situations:
  \begin{enumerate}
  \item We add an edge and its two endpoints. In that case, we are
    adding a new connected component which consists of a tree, which
    we know to be compatible.
  \item We add an edge $(u,v)$, and one of its two endpoints ($v$),
    the other one already being in $G$. We can extend $f_E$ by
    assigning any value to $f_E((u,v))$. We extend $f_V$ by setting
    $f_V(v) = f_E((u,v)) - f_V(u) \pmod{d}$. We then get that the
    augmented graph is additive, and hence compatible.
  \item We add an edge between two nodes of $G$. Let $f_V$ be a valid
    v-labeling of $(G,f_E)$. We extend $f_E$ to the new edge $(u,v)$
    by setting $f_E((u,v)) = f_V(u) + f_V(v) \pmod{d}$. This shows
    that the augmented graph is additive, and hence compatible.
  \end{enumerate}
\end{proof}

\begin{corollary}
  Given a graph $G$ with \emph{some} of its edges labeled in \Zd by
  a function $f$. If the subgraph of $G$ induced by the domain of $f$,
  labeled with $f$, is compatible, then there is an extension $f_E$ of
  $f$, such that $(G,f_E)$ is a compatible e-labeled graph.
\end{corollary}
\begin{proof}
  The result follows from Lemma~\ref{lem:extension}. We first decide
  whether the subgraph induced by $f$ is compatible. If it is, we add
  the remaining edges of $G$ one by one.
\end{proof}

\section{An efficient additivity test}
\label{sec:additivity}

Theorems~\ref{theo:compeqadd} and~\ref{theo:numbers} give a
theoretical characterization of additive e-labeled graphs. These
results are not practical per se, since they involve verifying certain
conditions on all the cycles of a graph. In this section we develop a
polynomial algorithm to test for additivity.

We tackle this problem by studying $ A_G$, the incidence matrix of
$G$. That is, $A_G \in \Z^{n\times m}$ such that
\begin{equation*}
[A_G]_{i,j} = 
\begin{cases}
  1 \quad \text{if vertex $v_i$ is incident with edge $e_j$,}\\
  0 \quad \text{if not}
\end{cases}
\end{equation*}

We use the Smith Normal Form (SNF) $S$ of $A_G$ together with the left
and right multipliers $U,V$. Here, $U \in \Z^{n \times n}, V \in
\Z^{m\times m}, S \in \Z^{n \times m}$ have the following properties:
\begin{itemize}
\item[-]  $U$ and $V$ are unimodular,
\item[-] $S$ is a diagonal matrix, $s_{i,i} | s_{i+1,i+1}$ for all $i$, and
\item[-] $A_G \, = \,  U S V$.
\end{itemize}

The authors of~\cite{minorsIncidence} show that the SNF 
$S$ of $A_G$ is
\begin{equation}  \label{eq:SNF}
  \begin{bmatrix}
    D & 0
  \end{bmatrix},
\end{equation}
where $D = \text{diag}(1,\ldots,1,\alpha)$ is a diagonal matrix with
$1$'s in every entry but the last one, which we call $\alpha$. This last
entry is $0$ if $G$ is bipartite (i.e. has no odd cycles) of $2$ if it
is not.

\medskip

\begin{definition}
  Let $G=(V,E)$ be a graph and let $C$ be any cycle of $G$. We
  associate a vector $\omega_C \in \Z^{|E|}$ with $C$. We index the
  coordinates of $\omega_C$ using the edges of $G$.
  
  Label the consecutive edges of $C$
  \begin{equation}
    \label{eq:edgesC}
    e_1,e_2,\ldots,e_{k-1},e_k,
  \end{equation}
  with $e_1$ any edge of the cycle.
  If $C$ is an even cycle, we adjoin $(-1)^i$ to $e_i$:
  \begin{equation}
    \label{eq:annotatedeven}
    e_1,-e_2,\ldots,(-1)^i e_i,\ldots,e_{k-1},-e_{k}.
  \end{equation}
  If $C$ is an odd cycle and $d$ is even, we adjoin $d/2$ to each edge:
  \begin{equation}
    \label{eq:annotatedodd}
    \frac{d}{2}\,e_1,\ldots,\frac{d}{2}\, e_i,\ldots,\frac{d}{2}\,e_{k}.
  \end{equation}
  
  Since $C$ need not be a simple cycle, some edges may appear more
  than once in~(\ref{eq:edgesC}). Let $e'_1,\ldots,e'_r$ be the
  \emph{distinct} edges of $C$. For each distinct edge $e'_i$, we
  define $\omega_{e'_i}$ to be the sum of the numbers adjoined to each
  appearance of $e'_i$ in~(\ref{eq:annotatedeven})
  or~(\ref{eq:annotatedodd}). For example, if an edge $e'_i$ appears
  twice, both times accompanied by a $1$, then the corresponding
  $\omega_{e'_i}$ is $2$. If one of the appearances has a $1$ and the other
  one a ($-1$), then $\omega_{e'_i}$ is $0$.

  Given a cycle $C$, we define $\omega_C$ as
  \begin{equation}
    \label{eq:omegac}
    (\omega_C)_{(u,v)} =
    \begin{cases}
      \omega_{(u,v)} & \text{if $(u,v)$ is in $C$,}\\
      0 & \text{otherwise}.
    \end{cases}
  \end{equation}

  Notice that in~(\ref{eq:annotatedeven}), the choice of $e_1$ may
  swap the $1$'s and the $-1$'s. This is not problematic, since it
  only changes $\omega_C$ into $-\omega_C$. The $\omega_C$, with $C$
  of even length, are in the kernel of the incidence matrix of $G$ and
  we only use them in that context.
\end{definition}

\begin{lemma}
  \label{lem:sumofcomponents}
  Let $C$ be a cycle of $G$. If the length of $C$ is even, then the
  sum of the coordinates of $\omega_{C}$ is $0$. If the length of $C$
  is odd, then the sum of the coordinates of $\omega_C \equiv d/2
  \pmod{d}$.
\end{lemma}
\begin{proof}
  Notice that we have the same number of edges accompanied by $1$ as
  the number of edges accompanied by $-1$. The sum of the coordinates
  of $\omega_{C}$ is the sum of all these $1$'s and $-1$'s, and
  is therefore $0$.

  The sum of the coordinates of $\omega_C$ is an odd integer multiple
  (i.e. the number of edges of $C$) of $d/2$.
\end{proof}

Let $(G,f_E)$ be an e-labeled graph, $\omega \in \Z^{|E|}$. We denote
\begin{equation}\label{eq:inner}
  \langle \omega,f_E \rangle := \sum_{(u,v) \in E} \omega_{(u,v)}
  f_E((u,v)).
\end{equation}

Let $\pi_d :\Z^{|E|} \to \Zd^{|E|}$ denote the projection
\begin{equation*}
  \pi_d(x)_{(u,v)} = r_d(x_{(u,v)}),
\end{equation*}
where $r_d$ is the remainder modulo $d$.  Finally, we denote by
$\mathcal{C}$ the set of even cycles in $G$.

The integer kernel of $A_G$ is computed in~\cite{rafael}, and is
shown to be the submodule spanned by $\{\omega_C, C\in\mathcal{C}\}$:
\begin{equation}
  \label{eq:kerZ}
  \ker_\Z(A_G) = \langle \omega_C, C\in\mathcal{C}\rangle.
\end{equation}
We prove a modular version of this result. Given $M \in \Z^{a\times
  b}$, we define $\ker_{\Zd}(M) = \{\mathbf{x} \in \Zd^b, M\mathbf{x}
\equiv 0 \pmod{d}\}$.

\begin{proposition}
  \label{prop:modulark}
  Let $G$ be a connected graph, and let $A_G$ be its incidence
  matrix. Then
  \begin{enumerate}
  \item If $d$ is odd or if $G$ has no odd cycles, then
    $\label{eq:kerZdodd} \ker_{\Zd}(A_G) = \pi_d(\ker_\Z(A_G))$.
  \item If $d$ is even and there is an odd cycle $C'$ in $G$, then
    \begin{equation*}
      \label{eq:kerZdeven}
      \ker_{\Zd} (A_G) = \pi_d(\ker_\Z(A_G)) \oplus
      \langle \pi_d(\omega_{C'}) \rangle.
    \end{equation*}
  \end{enumerate}
\end{proposition}
\begin{proof}
  In this proof, $\{z_1,\ldots,z_m\}$ denotes the canonical basis of
  $\Z^m$. That is, $(z_i)_i = 1$ and $(z_i)_j = 0$, for $j\neq
  i$. Analogously, $\{\pi_d(z_1),\ldots,\pi_d(z_m)\}$, denotes the
  canonical basis of $\Zd^m$.

  Let $S$ be the SNF of $A_G$, and $U$,$V$ such that $A_G = USV$, as
  described in \eqref{eq:SNF}. Equivalently, $U^{-1}A_G = SV$. Since
  $U$ and $V$ are both unimodular, they have integer inverses
  \emph{and} they have integer inverses modulo $d$. Therefore
  $\ker_\Z(A_G) = \ker_\Z(SV)$ and $\ker_\Zd(A_G) = \ker_\Zd(SV)$,
  implying that
  \begin{gather}
    \label{eq:kernelsAandS}
    \ker_\Z(A_G) = V^{-1}\ker_\Z(S)\\
    \label{eq:kernelsAandSd}
    \ker_\Zd(A_G) = \pi_d(V^{-1}\ker_\Zd(S))
  \end{gather}

  Let $\mathbf{x} = (x_1,\ldots,x_m) \in \ker_\Zd(S)$. That means that
  \begin{equation}
    \label{eq:Sx}
    S\mathbf{x} = 
    \begin{pmatrix}
      x_1\\
      \vdots\\
      \alpha x_n
    \end{pmatrix}
    \equiv 0 \pmod{d}
  \end{equation}

  If $\alpha = 0$ (i.e. $G$ has no odd cycles), equation~(\ref{eq:Sx})
  holds if and only if $x_1 = \cdots = x_{n-1} = 0$. That means that 
  \begin{equation*}
    \ker_\Zd(S) = \langle
    \pi_d(z_n),\ldots,\pi_d(z_m)\rangle \quad \text{and} \quad
    \ker_\Z(S) = \langle z_n,\ldots,z_m\rangle.
  \end{equation*}
  Therefore, we have $\ker_\Zd(S) = \pi_d(\ker_\Z(S))$, whence
  $\ker_\Zd(A_G) = \pi_d(\ker_\Z(A_G))$.
  
  If $\alpha = 2$ (i.e. $G$ has an odd cycle) and $d$ is odd,
  equation~(\ref{eq:Sx}) holds if and only if $x_1 = \cdots =
  x_n = 0$. That means that
  \begin{equation*}
    \ker_\Zd(S) = \langle \pi_d(z_{n+1}),\ldots,\pi_d(z_m)\rangle.
  \end{equation*}
  Once more, 
  \begin{equation*}
    \ker_\Z(S) = \langle z_{n+1},\ldots,z_m\rangle.
  \end{equation*}
  And again $\ker_\Zd(S) = \pi_d(\ker_\Z(S))$, implying $\ker_\Zd(A_G)
  = \pi_d(\ker_\Z(A_G))$.

  We now assume that $\alpha = 2$ and that $d$ is even. From
  equation~(\ref{eq:Sx}) we now deduce that
  \begin{gather}
    \label{eq:kerdS}
    \ker_\Zd(S) = \langle \pi_d(z_{n+1}),\ldots,\pi_d(z_m)\rangle \oplus \langle
    \frac{d}{2} \pi_d(z_n) \rangle,\\
    \label{eq:kerS}
    \ker_\Z(S) = \langle z_{n+1},\ldots,z_m\rangle.
  \end{gather}
  Notice that
  \begin{equation}
    \label{eq:smallSet}
    \langle \frac{d}{2} \pi_d(z_n) \rangle = \{ 0,\frac{d}{2}
    \pi_d(z_n)\}.
  \end{equation}
  Combining
  equations~(\ref{eq:kernelsAandS}),~(\ref{eq:kernelsAandSd}),~(\ref{eq:kerdS})
  and~(\ref{eq:kerS}), we have
  \begin{gather}
    \label{eq:kerdA}
    \ker_\Zd(A_G) = \langle
    \pi_d(V^{-1}\pi_d(z_{n+1})),\ldots,\pi_d(V^{-1}\pi_d(z_m))\rangle \oplus
    \langle \pi_d(V^{-1}\frac{d}{2} \pi_d(z_n)) \rangle.\\
    \label{eq:kerA}
    \ker_\Z(A_G) = \langle V^{-1}z_{n+1},\ldots,V^{-1}z_m \rangle.
  \end{gather}

  Let $C'$ be an odd cycle of $G$. Then $\pi_d(\omega_{C'})
  \in\ker_\Zd(A_G)$. To see why, recall that entry $e_j$ of
  $\omega_{C'}$ is $d/2$ times the number of occurrences of the edge
  $e_j$ in $C'$. For every vertex $v_i$ of the cycle, the number of
  edges that enter and leave it must be the same. That means that the
  $v_i$-th entry of $A_G\omega_{C'}$ has an even number times $d/2$
  (if vertex $v_i$ is in the cycle) or $0$. In both cases,
  $A_G\omega_{C'} \equiv 0 \pmod{d}$.

  Now, since $\pi_d(\omega_{C'}) \in \ker_\Zd(A_G)$, we must have
  \begin{equation}
    \label{eq:decompOmegaC}
    \pi_d(\omega_{C'}) = \sum_{l = n+1}^m \gamma_l \pi_d(V^{-1}z_l) +
    \varepsilon \pi_d(V^{-1}\frac{d}{2} \pi_d(z_n)),
  \end{equation}
  where $\varepsilon$ is $0$ or $1$ (see~(\ref{eq:smallSet})). The
  first summand consists of multiples of the projections of even
  cycles (see~(\ref{eq:kerZ})). That means that if we take the sum of
  the coordinates of both sides of equation~(\ref{eq:decompOmegaC}),
  we get $\varepsilon = 1$ (see Lemmas~\ref{lem:sumofcomponents}
  and~\ref{lem:sumofcomponents}.) If we set
  \begin{equation*}
    \gamma = \sum_{l = n+1}^m \gamma_l \pi_d(V^{-1}z_l),
  \end{equation*}
  we can write
  \begin{equation}
    \label{eq:rewriteomegaC}
    \pi_d(V^{-1}\frac{d}{2} \pi_d(z_n)) = \gamma - \pi_d(\omega_{C'})
  \end{equation}

  Now, take any $\mathbf{x} \in \ker_\Zd(A_G)$. We have that
  \begin{equation*}
    \mathbf{x} = \sum_{l = n+1}^m \beta_l \pi_d(V^{-1}z_l) +
    \beta \pi_d(V^{-1}\frac{d}{2} \pi_d(z_n)).
  \end{equation*}
  Plugging in equation~(\ref{eq:rewriteomegaC}) we get 
  \begin{equation*}
    \mathbf{x} = \sum_{l = n+1}^m \beta_l
    \pi_d(V^{-1}z_l) +
    \beta (\gamma - \pi_d(\omega_{C'})).
  \end{equation*}
  If we set $\tilde{\beta_l} = \beta_l + \gamma_l$, we have
  \begin{equation*}
    \mathbf{x} = \sum_{l = n+1}^m
    \tilde{\beta_l} \pi_d(V^{-1}z_l) + (-\beta)\pi_d(\omega_{C'}),
  \end{equation*}
  which shows that
  \begin{equation*}
    \ker_{\Zd} (A_G) = \pi_d(\ker_\Z(A_G)) \oplus
    \langle \pi_d(\omega_{C'}) \rangle.
  \end{equation*}
\end{proof}

The results we have discussed allow us to obtain the following
\begin{theorem}
  \label{theo:kernels}
  Let $(G,f_E)$ be an e-labeled connected graph. Let $A_G$ be the
  incidence matrix of $G$. The following statements are equivalent:
  \begin{enumerate}
  \item $(G,f_E)$ is a compatible e-labeled graph.
  \item $\langle \pi_d(\omega_C),f_E \rangle \equiv 0 \pmod{d}$, for every
    cycle $C$ of $G$.
  \item $\langle \omega,f_E \rangle \equiv 0 \pmod{d}$, for
    all $\omega \in \ker_{\Zd}(A_G)$.
  \item If $d$ is odd or $G$ has no odd cycles, $\langle \omega,f_E
    \rangle \equiv 0 \pmod{d}$, for all $\omega$ belonging to the
    projection of a finite set of generators of $\ker_{\Z}(A_G)$. If
    $d$ is even and has an odd cycle, $\langle \omega,f_E \rangle
    \equiv 0 \pmod{d}$, for all $\omega$ belonging to a finite set of
    generators of $\ker_{\Z}(A_G)$ and for $\omega_C$, for some odd
    cycle $C$.
  \item $(G,f_E)$ is an additive e-labeled graph.
  \end{enumerate}
\end{theorem}
\begin{proof}
  Clause~i) is equivalent to clause~v) by
  Theorem~\ref{theo:compeqadd}.  Clause~ii) is a restatement of
  clause~i) using a different notation. Clauses~ii) and~iii) are
  equivalent by Proposition~\ref{prop:modulark}. Clauses~iii) and~iv)
  also follow from that proposition: the finite sets described in
  clause~iv) were shown to be generators of $\ker_{\Zd}(A_G)$.
\end{proof}

The equivalence of clauses~iv) and~ v) in the above Theorem, provides
the following complexity result.

\begin{theorem}\label{th:complexity}
  Let $(G,f_E)$ be an e-labeled connected graph. The additivity of
  $(G,f_E)$ can be tested in time polynomial in the size of the
  graph. Furthermore, we can obtain all its valid v-labelings in
  polynomial time too.
\end{theorem}
\begin{proof}
  We compute the Smith Normal Form (SNF) $S$ of $A_G$ described in the
  proof of Proposition~\ref{prop:modulark}, together with the left and
  right multipliers $U$ and $V$. This computation can be carried out
  using the polynomial algorithm presented in~\cite{polySmith},
  modified to work with rectangular matrices in the way the authors of
  that paper suggest.

  We saw in Proposition~\ref{prop:modulark} that we can obtain
  generators of $\ker_{\Zd}(A_G)$ from the columns of $V^{-1}$. If $G$
  has no odd cycles (i.e. $\alpha = 0$), we use the last $m-n+1$
  columns. If $\alpha = 2$ and $d$ is odd, we use the last $m-n$
  columns. If $\alpha = 2$ and $d$ is even, we use the last $m-n$
  columns and $d/2$ times its $n$-th column. To check the additivity
  of $(G,f_E)$, we just need to verify that these generators satisfy
  the conditions stated in clause iv) of Theorem~\ref{theo:kernels}.

  Once we know that $(G,f_E)$ is additive, we can efficiently obtain
  all its valid v-labelings. We must first know whether $G$ has an odd
  cycle or not. This can be read directly off the SNF $S$ of $A_G$:
  $G$ has an odd cycle if and only if the diagonal of $S$ contains a
  $2$. Having no odd cycles is classically known to be equivalent to
  $G$ being bipartite (cf. for instance ~\cite{harary}, p. 18), and
  can be checked in time $O(n+m)$. We can also obtain an odd cycle in
  $G$ as a byproduct of this check.

  If $G$ has no odd cycles, we can assign any of the $d$ possible
  labels to an arbitrary vertex, and then propagate the label to the
  rest of the graph using breadth-first search (BFS). If $G$ does have
  odd cycles, choose one of them and call it $C$. Choose a vertex
  $v_1$ in $C$. Formula~(\ref{eq:condodd}) shows which label (or
  labels, if $d$ is even) we can assign to $v_1$ in order to obtain
  valid v-labelings of $(G,f_E)$. We then propagate the label of $v_1$
  to the rest of the graph using BFS.
\end{proof}

\begin{remark}
  Given a graph $G$, consider the \emph{cycle space} of $G$
  (\cite{harary}). It is the $\Z_2$-vector space generated by the
  fundamental cycles of $G$. That is, the cycles obtained when adding
  an edge of $G$ to a spanning tree.

  One might be tempted to think that checking the compatibility
  conditions on these generators suffices to verify the compatibility
  of a graph with labels in \Zd for any $d$, as in the case $d=2$. However, consider for
instance the graph in Figure~\ref{fig:cycleSpace}, in which we marked the spanning tree
with edges $\{e_{14}, e_{23}, e_{24}\}$:
  \begin{figure}[ht]
    \centering
    \includegraphics{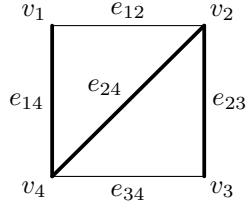}
    \caption{A spanning tree of a graph.}
    \label{fig:cycleSpace}
  \end{figure}
The sum of the two fundamental triangle cycles $C_1, C_2$ (represented by their $0,1$ vectors) 
equals the square cycle $C$ only when $d=2$.
This
  situation is depicted informally in Figure~\ref{fig:sum}.
  \begin{figure}[ht]
    \centering
    \includegraphics{additive17.eps}
    \caption{Adding two odd cycles to obtain an even one.}
    \label{fig:sum}
  \end{figure}
However, if $d$ is odd we do not impose any conditions on $C_1$ and $C_2$, and so this cannot insure
the even cycle condition we need to check. When $d\neq 2$ is even, we  get $\frac d 2$ times the
even cycle condition, which again is not sufficient to insure additivity. Consider for instance 
the labeling $f_(e_{12}) = f_(e_{24}) =f_(e_{34}) = 1, \, f_(e_{14}) =f_(e_{23)} = 0$ and $d=4$. The odd cycle property is
verified for $C_1, C_2$ but the labeling is not additive.
\end{remark}

\section{Multiplicative version}

In the previous sections, we used labelings that assigned integers
modulo $d$ to the edges and vertices of a graph. But actually,
everything we wrote is also valid if the labels belong to any finite
cyclic group, via the isomorphism with \Zd. In particular, we can use
labelings in \Gd, the $d$-th roots of unity. In this case, the
isomorphism between $\Zd$ and $\Gd$ is given by
\begin{equation}
  \label{eq:isomorphism}
  k \mapsto e^{2\pi ik/d}.
\end{equation}
This alternate formulation is useful because it links our problem with
the theory of toric ideals. As a general text on this subject, we
refer the reader to~\cite{convexpoly}.

Let us state this equivalent version. Let $G=(V,E)$ be a connected graph and
$d$ an integer greater than $1$. Let $n = |V|$ and $m = |E|$. Let
$v_1,\ldots,v_n$ be the vertices of $G$ and let $e_1,\ldots,e_m$ be
its edges. We work with complex variables $x_{v_i}$ for each $v_i \in
V$, and $y_{e_i}$ for each $e_i \in E$. The value of $x_{v_i}$
corresponds to the label of vertex $v_i$, and the value of $y_{e_i}$
corresponds to the label of edge $e_i$. We can restate
Problem~\ref{prob:additive} in this multiplicative setting:

\begin{problem}\label{prob:multiplicative}
  For which $\mathbf{y} \in \Gd^m$ are there $\mathbf{x} \in \Gd^n$
  such that
  \begin{equation}
    \label{eq:multcond}
    y_{e_i} = x_{u_i}x_{v_i},
  \end{equation}
   holds for every edge $e_i = (u_i,v_i) \in E$?
\end{problem}

According to a classic result for toric parametrizations, given a
vector $\mathbf{y} \in (\C^*)^m$ of complex nonzero numbers, there is
an $\mathbf{x} \in (\C^*)^n$ satisfying~(\ref{eq:multcond}) if and
only if

\begin{equation}
  y^{\mathbf{u}} = y_1^{u_1}\cdots y_m^{u_m} = 1,
\end{equation}
for every $\mathbf{u} = (u_1,\ldots,u_m) \in
\ker_{\Z}(A_G)$. Furthermore, when these conditions are satisfied,
the number of such solutions is
\begin{equation} \label{eq:gcd}
  g = \gcd(\{\text{maximal minors of $A_G$} \}),
\end{equation}
provided that $g\neq 0$, in which case there are infinitely many
solutions. We deduce from~\eqref{eq:SNF} that $g = 2$ or $0$,
depending on whether $G$ has an odd cycle or not, respectively. It was
this result which prompted us to study the incidence matrix of $G$ in
connection with Problem~\ref{prob:additive}.

We now state a modular version of the toric result.  We impose the
additional restriction that
\begin{equation}
  \label{eq:aredthroots}
  x_{v_i}^d = 1,
\end{equation}
for all $v_i \in V$.  This condition, together
with~(\ref{eq:multcond}), implies that the $y_{e_i}$ are also in \Gd.

\begin{theorem}
  \label{theo:modulartoric}
  Let $G = (V,E)$ be a connected graph. Given $\mathbf{y} \in \Gd^m$,
  there exists $\mathbf{x} \in \Gd^n$ satisfying (\ref{eq:multcond}) if
  and only if
  \begin{equation}
    y^{\mathbf{u}} = 1,
  \end{equation}
  for every $\mathbf{u} \in \ker_{\Zd}(A_G)$. If $g$ is $0$, there are
  $d$ solutions to~(\ref{eq:multcond}) and~(\ref{eq:aredthroots})
  simultaneously. If $g$ is $2$ and $d$ is even, there are two
  solutions. Otherwise, there is a unique solution.
\end{theorem}


The result can be translated from Theorem~\ref{theo:numbers}.
Alternatively, we could prove that given $\mathbf{y} \in \Gd^m$, there
are as many solutions $\mathbf{x} \in \Gd^n$ as stated using the
knowledge of $g$ in \eqref{eq:gcd}, by checking how many of the
complex solutions $\mathbf{x} \in (\C^*)^n$ consist of $d$-th roots of
unity.

\bibliography{additive}

\bibliographystyle{plain}

\end{document}